\documentclass{article} \usepackage{amssymb}\usepackage{amscd}
\usepackage[utf8]{inputenc}
\usepackage{amsfonts}

\usepackage{amsmath}

\usepackage{color}
\usepackage{amstext}
\newtheorem{theorem}{Theorem}
\newtheorem{lemma}{Lemma}\newtheorem{claim}{Claim}

\newtheorem{problem}{Problem}
 \newtheorem{remark}{Remark}\newtheorem{proposition}{Proposition}
\newtheorem{corollary}{Corollary} \newcommand{\La}{\Lambda}
\newcommand{\R}{{\mathbb R}}  \newcommand{\Z}{{\mathbb Z}} \newcommand{\N}{{\mathbb N}}

\newcommand{\la}{\lambda}
\newcommand{\de}{\delta}
\newcommand{\dist}{{\rm dist}}
\newcommand{\Cc}{{\mathbb C}}

\newcommand{\pw}{$PW_S^1$\ }
\newcommand{\ext}{Ext$(PW_S^1)$}
\newcommand{\expo}{Exp$(PW_S^1)$}

\DeclareMathOperator*{\Res}{Res}

\begin{document}
\title{On geometry of the unit ball of Paley-Wiener space over two symmetric intervals}
\author{A. Ulanovskii \and I. Zlotnikov}\date{}
\maketitle
\begin{abstract}Let $PW_S^1$ be the space of integrable functions on $\R$ whose Fourier transform vanishes outside $S$, where $S = [-\sigma,-\rho]\cup[\rho,\sigma]$, $0<\rho<\sigma$. In the case $\rho>\sigma/2$, we present a complete description of the set of extreme and the set of exposed points of the unit ball of \pw\!. The structure of these sets becomes more complicated when $\rho<\sigma/2$. 
\end{abstract}

\section{Introduction}

Given a Banach space $X$, denote by
$$\mbox{ball}(X) := \{f\in X:\|f\|\leq 1\}$$
the closed unit ball of $X$. As usual, an element $f$ from ball$(X)$ is called {\it extreme} if it is not an interior point of any line segment contained in ball$(X).$
An element $f$
in ball$(X)$  is an {\it exposed point} of ball$(X)$, if there exists a functional $\phi\in X^\ast$ such that $\|\phi\|= 1$ and the set $\{g\in X : \phi(g)=1\}$ consists of one
element, $f$. It is easy to check that every exposed point is extreme.

The classical theorem of K.~de~Leeuw and W.~Rudin \cite{dLR} states that the extreme points of the unit ball of the Hardy space $H^1$ on the unit disk are
precisely the outer functions $f\in H^1$ with $\|f\|_1 = 1$, see also \cite{Ga}, Chapter IV and \cite{H}, Chapter 9 for alternative presentations.
See also \cite{D} for an extension to `punctured'  Hardy  spaces.  On the other hand, no description of the exposed points of ball$(H^1)$ is known, though this set has been studied by
a number of authors, see e.g. \cite{p} and the literature therein. 
 The extreme and exposed points of unit ball have been determined for certain spaces of polynomials equipped with $L^1$-norm (\cite{D_PW},  \cite{D_lacun}), entire functions of exponential type  equipped with $L^1$-norm (\cite{D_PW}), and in some other spaces. See also the  list of references in the above papers.

Given a function $f\in L^1(\R)$, denote by
$$\mbox{Sp}(f):=\overline{\{x\in\R: \hat f(x)\ne0\}}$$   the (closed) spectrum of $f$. Here $\hat f$ is the Fourier transform of $f,$
$$
\hat f(t)=\int\limits_\R e^{-2\pi i t x}f(x)\,dx.
$$Below we will use the  well-known fact that the spectrum is  defined for any tempered distribution.

Let $S$ be  a finite or infinite union of disjoint closed intervals. Denote by  \pw the Paley-Wiener space
$$
PW_S^1:=\{f\in L^1(\R):\mbox{ Sp}(f)\subset S\}
$$equipped with the $L^1$-norm $\|\cdot\|_1$. By \ext\, and \expo\! we denote the set of extreme and exposed points of ball(\pw), respectively.

Observe that when  $S=[-\sigma,\sigma],\sigma>0, $ is a single interval, the space \pw\, consists of  entire functions $f\in L^1(\R)$ satisfying the inequality
$$
|f(x+iy)|\leq Ce^{2\pi\sigma |y|}, \quad x,y\in\R,
$$where $C=C(f)$ is a constant. This implies that  $f$ is an {\it entire function of exponential type}  $\leq 2 \pi \sigma$. We will denote by Type$(f)$ the exponential type of $f$, see definition in Sec. 2.1.


We are interested in the following

\begin{problem}[K. Dyakonov \cite{D_lacun}] Describe  the sets {\rm \ext} and {\rm \expo}.\end{problem}

 Denote by Hol$(\Cc)$ the set of all entire functions. Assume  $f\in$Hol$(\Cc)$. Let us introduce two sets which will play an important role in what follows:
 \medskip

(i) Denote by $\La(f)\subset\Cc\setminus\R$   the (possibly empty)  set of all  points $\lambda=a+ib, a,b\in\R,b\ne0,$ such that
$f(\lambda)=f(\bar\lambda)=0$, where $\bar\lambda=a-ib.$

(ii) Denote by  $\Omega(f)\subset\R$   the (possibly empty)  {\it multiset}  of all  points $x\in \R$ such that $f(x)=f'(x)=...=f^{(2n(x)-1)}(x)=0$, where $n(x)\geq1$ denotes the integer such that $f$ has  zero of multiplicity $ 2n(x)$ or $2n(x)+1$ at $x$.  Every point $x\in\Omega(f)$ is {\it counted with multiplicity $2n(x)$}.

\medskip

If $S=[-\sigma,\sigma]$ is a single interval, then  Problem 1 has the following solution:

\begin{theorem}[\cite{D_PW}]\label{t0} {\rm (A)}  The set {\rm Ext}$(PW_{[-\sigma,\sigma]}^1),\sigma>0,$ consists precisely of the functions $f$ satisfying the conditions:
\begin{equation}\label{c0}\|f\|_1=1;\end{equation}
\begin{equation}\label{d1}
\mbox{ at least one of the points  }\{\pm \sigma\}\in \mbox{{\rm Sp}}(f);
\end{equation}
\begin{equation}\label{2}
\La(f)=\emptyset.
\end{equation}

{\rm (B)} The set {\rm Exp}$(PW_{[-\sigma,\sigma]}^1),\sigma>0,$ consists precisely of the functions $f$ from {\rm Ext}$(PW_{[-\sigma,\sigma]}^1)$ satisfying the conditions:
\begin{equation}\label{1exp}
\Omega(f)=\emptyset;\end{equation}
\begin{equation}\label{exp2}
\int_\R |f(x)|w(x)dx = \infty,  \forall w\in\mbox{{\rm Hol}}(\Cc),\mbox{{\rm Type}}(w)=0, w|_\R\geq0, w\ne const.
\end{equation} \end{theorem}

Condition (\ref{2}) means that no function $f\in$\ext\, may have symmetric (with respect to $x$-axis) zeros. Condition  (\ref{1exp}) means that no function $f\in$\expo\, may have a real zero of order $\geq 2$.
We note also that condition (\ref{exp2}) means that $fw\notin L^1(\R)$, for every entire real (on $\R$) nonconstant function $w$ of zero type.

In this paper we consider the  spectra $S$ which consist of two intervals. For simplicity of presentation, we will always assume that the intervals have the same length and are symmetric:
$$
S:=[-\sigma,-\rho]\cup[\rho,\sigma]=[-\sigma,\sigma]\setminus(-\rho,\rho),\quad 0<\rho<\sigma,
$$though some of the results remain valid for more general situation, see Theorem~\ref{tlast} in Sec. 2. We will  say that  $S$ has the gap $(-\rho,\rho)$.

Clearly, condition (\ref{c0}) is necessary for $f$  to be an extreme point. Observe also that since  \pw$\subset PW_{[-\sigma,\sigma]}^1$, the inclusion holds true:
\begin{equation}\label{ex_inclusion}
    \mbox{{\rm \pw}}\!\cap\mbox{{\rm Ext}}(PW^1_{[-\sigma,\sigma]})\subset\mbox{{\rm \ext}}.
\end{equation} One may ask if the above inclusion is proper, i.e. (see part (A) of Theorem \ref{t0}) if the set \ext\, may contain functions $f$ that do not satisfy  (\ref{d1}) and (\ref{2}).  We will see that this is the case. A similar remark is true for the set of exposed points. Moreover, the results below show that the structure  of both \ext\! and \expo\!
depends on the size of the gap $(-\rho,\rho)$.

Observe that  every function $f\in PW_S^1$ admits a unique representation
$$f(z)=f_-(z)+f_+(z),\quad \mbox{Sp}(f_-)\subset[-\sigma,-\rho],\ \mbox{Sp}(f_+)\subset[\rho,\sigma].$$
It is easy to check that $f_-\in PW_{[-\sigma,-\rho]}^1$ and $f_+\in PW_{[\rho,\sigma]}^1$, though we do not use this fact below.


If $\rho>\sigma/2$ (`long gap'), then the description of \ext\, and \expo\! is somewhat similar to the one in Theorem \ref{t0}:

\begin{theorem}[Long gap]\label{t1}
Assume that  $\sigma/2<\rho<\sigma.$

{\rm (A)} The set  {\rm \ext} consists precisely of the functions $f$ satisfying {\rm(\ref{c0})} and the conditions   \begin{equation}\label{sp_cond}\mbox{at least one of the points }\, \{\pm \sigma\},\{\pm\rho\}\mbox{ belongs to } \mbox{{\rm Sp}}(f); \end{equation}\begin{equation}\label{la_cond}
\La(f_-)\cap \La(f_+)=\emptyset.
\end{equation}

{\rm (B)} The set  {\rm \expo} consists precisely of the functions $f$ from {\rm \ext} satisfying  {\rm (\ref{exp2})} and    \begin{equation}\label{exp1}
\Omega(f_-)\cap\Omega(f_+)=\emptyset.\end{equation}
\end{theorem}

One may check that condition (\ref{la_cond}) does not imply condition (\ref{2}). Similarly, condition (\ref{exp1}) does not imply condition (\ref{1exp}).

Observe that condition (\ref{exp1}) implies \begin{equation}\label{omega}\mbox{{\rm There is no }} x\in\R \ \mbox{{\rm  such that }} f_-(x)=f'_-(x)=f_+(x)=f'_+(x)=0.\end{equation}

We also show that  the set of extreme points \ext \, is `large':

 \begin{theorem}[Long gap]\label{tt0} Assume that
 $\sigma/2<\rho<\sigma$.

 {\rm (i)} If $f\in$\pw\! and  $\|f\|_1 = 1$,  then  there exist extreme points $f_1$ and $f_2$
  such that $f_1+f_2=2f$.

  {\rm (ii)} The set {\rm \ext}\! is dense on the unit sphere $\{g\in$\pw\!$:\|g\|_1=1\}. $
\end{theorem}

 An analogue of Theorem \ref{tt0} is true for the Paley-Wiener space $PW^1_{[-\sigma,\sigma]}$ (see Corollaries 1 and 2 in \cite{d1}).
Observe also that part (i)  is true  for the Hardy space $H^1$ (see Theorem 2 in \cite{dLR}).

The proofs of necessity of the conditions  (\ref{sp_cond}) and (\ref{la_cond}) in part (A) and the conditions (\ref{exp2}), (\ref{exp1})
in part (B) of Theorem~\ref{t1} are easy. Moreover, these conditions remain necessary for every $\rho,0<\rho<\sigma,$ see Proposition \ref{p2} below.

The proof of sufficiency is more involved. The main step  (see Corollary \ref{c00}) is to show   that
 the set  $\La(f)\cup\Omega(f)$ is `small' with respect to the size of gap in the following sense: For every function $f\in PW_S^1$,  the  exponential system  with frequencies in $\La(f)\cup\Omega(f)$ (see definition below) is not complete in $L^2$ on some proper  subinterval of $(-\rho,\rho)$.
Moreover, the sufficiency part of Theorem~\ref{t1} remains true  (for every spectra $S$) for the functions $f\in$\pw\!  satisfying (\ref{d1}) and  the incompleteness property above, see Corollary~\ref{carb} in Sec.~2.
However, the incompleteness property   fails  in general to be true  for $\rho<\sigma/2$. Indeed, Theorem \ref{t1} does not hold  in this case even for functions satisfying (\ref{d1}):

\begin{theorem}[Short gap]\label{t2}
Assume that $ 0<\rho<\sigma/2.$

{\rm (i)} There is a function  $f$ from  $PW_S^1$ satisfying {\rm(\ref{c0})},    {\rm(\ref{la_cond})}, and
{\rm (\ref{d1})} such that $f\not\in\,${\rm \ext}.

{\rm (ii)} There is a function  $f$ from  {\rm \ext\!} satisfying  {\rm (\ref{exp2})}, {\rm (\ref{exp1})}, and {\rm (\ref{d1})} such that $f\not\in\,${\rm \expo}.
\end{theorem}

 Theorems \ref{t1}  and \ref{t2}  show that the value $\rho=\sigma/2$ is  `critical' for the structure  of the sets \ext\, and \expo.
  For the functions $f$ satisfying (\ref{d1}),  the reason for that is that the sets $\La(f)$ and $\Omega(f)$ can be `large' with respect to the size of gap. For functions $f$ that do not satisfy (\ref{d1}), there is another simple reason for why the sufficiency part of Theorem \ref{t1} ceases to be true in the case of short gap, see Sec. 5.

   However,   the size of the sets $\La(f)$ and $\Omega(f)$ alone does  not determine  if $f$ is extreme (exposed) or not:

\begin{theorem}[Short gap]\label{t3}

{\rm (i)} Assume that $0<\rho<\sigma/4.$ Then there exist   functions  $f_1,f_2\in PW_S^1$ satisfying {\rm(\ref{c0})},  {\rm(\ref{la_cond})}, and {\rm (\ref{d1})} such that~$$\La(f_1)=\La(f_2), \ f_1\in\mbox{{\rm \ext}}, \ f_2\not\in\mbox{\rm {\ext}}.$$

{\rm (ii)} Assume that $0<\rho<\sigma/5.$ Then there exist   functions  $f_1,f_2\in${\rm \ext}  satisfying {\rm(\ref{exp2})}, {\rm(\ref{exp1})}, and {\rm (\ref{d1})}  such that$$\Omega(f_1)=\Omega(f_2), \ f_1\in\mbox{{\rm \expo}}, \ f_2\not\in\mbox{\rm {\expo}}.$$
\end{theorem}




\section{Long gap: Proof of Theorem \ref{t1}}

\subsection{Auxiliary results}

Recall that an entire function $f$ is  of exponential type $d\geq0$ if for every positive $\varepsilon$ there is a constant $C_\varepsilon$ such that
$$
|f(z)|\leq C_\varepsilon e^{(d+\varepsilon)|z|},\quad z\in\Cc,
$$and
$$
|f(z)|>e^{(d-\varepsilon)|z|},\quad \mbox{for some sequence } z=z_n\in\Cc, |z_n|\to\infty.
$$We will write $d=\,$Type$(f)$.

We denote by Cart the {\it Cartwright class} which consists of all entire functions $f$ satisfying Type$(f)<\infty$ and
$$
\int_\R\frac{\log^+|f(x)|}{1+x^2}\,dx<\infty.
$$

For every $f\in \,$Cart  there exist numbers $d_\pm$ such that the following relations
$$
\log|f(z)|=d_- y+o(|z|),\quad z=x+iy, y\leq0;
$$$$
\log|f(z)|=d_+y +o(|z|),\quad z=x+iy, y\geq0
$$hold outside a `small' set, see  \cite{levin}, Lec. 16, Theorem 2.
Set$$
 \mbox{{\rm Type}}_\pm(f):=\limsup_{y\to\pm\infty}\frac{\log|f(iy)|}{|y|}.
$$Then  the limits above exist provided $y$ lies outside a `small' set, and we have Type$_\pm(f)=d_\pm,$ see \cite{levin}, Lec. 16.

The following lemma is well-known:

\begin{lemma}\label{ko} Assume  $f, g \in\,${\rm Cart}. Then

{\rm(i)} {\rm Type}$(f)= \max\{${\rm Type}$_-(f),${\rm Type}$_+(f)\}$;

{\rm(ii)} {\rm Type}$_\pm(fg) =${\rm Type}$_\pm(f)+${\rm Type}$_\pm(g)$;

{\rm(iii)} If  $g/f\in${\rm Hol}$(\Cc)$, then $g/f\in${\rm Cart} and $$\mbox{{\rm Type}}_\pm(g/f) =\mbox{{\rm Type}}_\pm(g)-\mbox{{\rm Type}}_\pm(f).$$

\end{lemma}

This lemma was also used in \cite{D_PW}. See \cite{hj}, pp. 192--194 for a proof of (i) and (iii), and \cite{levin}, Lec. 16, Theorem 4, for a proof of (ii).

Observe also that by the Phragm\'{e}n-Lindel\"{o}f theorem (\cite{levin}, Lec. 6, Theorem~3), every function $f$ of exponential type $d$ bounded on $\R$ satisfies the inequality \begin{equation}\label{in}|f(x+iy)|\leq Ce^{d|y|},\quad x,y\in\R, \end{equation}where the constant $d$ cannot be replaced with a smaller constant.

Given an entire function $f$, denote by $Z(f)$ the  sequence (multiset) of its zeros (counting multiplicities). Write $Z(f)=Z_r(f)\cup Z_c(f)$, where $Z_r(f)$ is the multiset of all real zeros and $Z_c(f)$ is the multiset of purely complex zeros of $f$.

 \begin{lemma}\label{l001}
Suppose functions $f,g$ from $PW^1_S$ are such that the ratio $h=g/f$ is  nonconstant and  real  on $\R$.

{\rm (i)}  Then the set of purely complex poles of $h$ is  a subset of $\La(f)$;

{\rm (ii)} If $h$ is bounded on $\R$, then $Z_r(f)\subset Z_r(g)$;

{\rm (iii)} If $h$ is an entire function of zero type bounded on $\R$, then $h=const.$

\end{lemma}

(i) Indeed, since  $h$ is real on $\R$, if a point $\lambda\in \Cc\setminus\R$ is a pole of $h$ then so is $\bar\lambda$. Therefore,
$\lambda\in\La(f).$

(ii) This is obvious.

(iii) This follows from (\ref{in}).


\medskip



Given a complex sequence $\Gamma\subset\Cc$, we denote by
$$
E(\Gamma):=\{t^n e^{2\pi i \gamma t}, \gamma\in\Gamma, n=0,...,m(\gamma) - 1\}
$$the corresponding exponential system, where $m(\gamma)\geq 1$ is the multiplicity of element $\gamma\in\Gamma$.
We denote by $R(\Gamma )$  the completeness radius of $E_\Gamma$,
$$R(\Gamma):=\sup\{a\geq0: E(\Gamma) \mbox{ is complete in } L^2(-a,a)\}.$$

We will use the following

\begin{lemma}[Extreme and exposed point criterion]\label{h_lemma}
{\rm (i)} Suppose that $f \in PW^1_S$ and $\|f\|_{1} = 1$.
Then $f\in${\rm Ext}$(PW^1_S)$ if and only if there is no real, nonconstant, bounded  function $h$ on $\R$ satisfying $fh \in PW^1_{S}$.

{\rm (ii)} Suppose that $f\in${\rm \ext}. Then $f\in${\rm \expo} if and only if there is no  nonnegative, measurable, nonconstant function $h$ on $\R$ satisfying $fh\in PW_S^1$.
\end{lemma}

We refer the reader to Lemma 1 in \cite{D_PW} for a proof of this lemma, see also \cite{G}, Chapter V.

 Lemma \ref{h_lemma} reduces Problem 1 to the following complex-analytic problems:

 \medskip
 (i) Determine for which functions $f\in$\pw\! there is no function $g\in$\pw\! such that the ratio $h:=g/f$ is  nonconstant, real and bounded on $\R$.

 (ii)  Determine for which functions $f\in$\ext\, there is no function $g\in$\pw\! such that the ratio $h:=g/f$ is  nonconstant and nonnegative on $\R$.

\medskip

The following proposition and its immediate corollary  play an essential role in solving the above problems for the case of long gap.

\begin{proposition}\label{p1}
Let $0<\rho<\sigma$ and $f\in$\pw\!.
Then $R(\La(f)\cup\Omega(f))\leq \sigma- \rho.$
\end{proposition}


\begin{corollary}[Long gap]\label{c00} Let $\sigma/2<\rho<\sigma$ and $f\in$\pw\!. Then  $R(\La(f)\cup\Omega(f))<\rho.$
\end{corollary}

The proof of Proposition \ref{p1} below uses  the classical Beurling-Malliavin theorem on completeness radius and a recent result on density of sign changes of real measures with spectral gap at the origin by M. Mitkovski and A. Poltoratski. Observe that  Proposition \ref{p1} remains valid for the larger class of entire functions that are tempered distributions whose spectrum lies in $S$.

\subsection{Beurling-Malliavin completeness theorem}
We will need several definitions.

 Given a point $\alpha\in\Cc$,  consider the mapping $T(\alpha)= \beta, \beta\in\R,$ where
$$
\frac{1}{\beta}=\frac{1}{2}\left(\frac{1}{\alpha}+\frac{1}{\bar \alpha}\right).
$$Clearly, $T(\alpha)=\alpha$, $\alpha\in\R,$ and $T(\alpha)=\infty$, Re$\,\alpha= 0$.

{\it Definition 1}.  A sequence of disjoint intervals
$\{I_n\}$ on the real line is called long (in the sense of Beurling and Malliavin)
if
$$
\sum_n \frac{|I_n|^2}{1+\dist^2(0,I_n)}=\infty,
$$where $|I_n|$ stands for the length of $I_n$. If the sum is finite, we call $\{I_n\}$ short.

{\it Definition 2}. Following \cite{BM_trick}, we say that a  sequence $\Gamma\subset\R$ is $d$-regular if for every
$\epsilon > 0$ any sequence of disjoint intervals $\{I_n\}$ that satisfies
$$
\left|\frac{\#(\Gamma\cap I_n)}{|I_n|}- d\right|\geq\epsilon,\quad \mbox{ for all } n,$$ is short.

There is a number of  slightly different definitions of $d$-regularity, see e.g. Remark 4.6 in \cite{{P_sp_gap}}. The definition above was used e.g. in \cite{MP}.

A complex sequence $\Gamma\subset\Cc\setminus (i\R)$ is called $d$-regular if the real sequence $T(\Gamma)$ is $d$-regular.

{\it Definition 3}. (i)
The interior Beurling-Malliavin density of a sequence $\Gamma\subset\R$ is defined as
$$D_{BM}(\Gamma) := \sup\{ d: \exists  d- \mbox{regular subsequence } \Gamma' \subset\Gamma\}.
$$

(ii) The exterior Beurling-Malliavin density can be defined as
$$D^{BM}(\Gamma):= \inf\{d: \exists d-\mbox{regular sequence } \Gamma'\supset\Gamma\}.$$

(iii) It will be more convenient for us to use an equivalent definition:
$$D^{BM}(\Gamma):= \sup\{ d: \exists \, \mbox{long sequence of intervals } I_n \mbox{ with }  \#(\Gamma\cap I_n)\geq d |I_n|,\forall n\}.$$

As an example, one may check that $D^{BM}(\Z)=D_{BM}(\Z)=1.$

\begin{theorem}[Beurling-Malliavin, \cite{BM_trick}]\label{bm} Let $\Gamma\subset\Cc$ be a sequence satisfying $$\sum_{\gamma\in\Gamma}\left|{\rm Im}\frac{1}{\gamma}\right|<\infty.$$ Then
\begin{equation}\label{R_DBM}
    R(\Gamma)=D^{BM}(T(\Gamma))/2.
\end{equation}
\end{theorem}

As a classical example, we recall that the trigonometric system $E(\Z)=\{e^{2\pi i nt},n\in\Z\}$ is complete in $L^2(-\rho,\rho)$ for every $\rho\leq1/2$, and is not complete when $\rho>1/2$, so that $R(\Z)=D^{BM}(\Z)/2=1/2.$

\subsection{Sign changes of real functions with spectral gap}
We say that a  function $f\in L^1(\R)$ (or a finite measure $\mu$) has a spectral gap at the origin  if there exists $r>0$ such that the Fourier transform of $f$ (of $\mu$) vanishes for $|t|\leq r$. Clearly,   every $f\in$\pw has spectral gap $[-\rho,\rho]$.

It is well-known that {\it real} functions and measures with spectral gap at the origin must have many oscillations (sign changes), see \cite{en} for definition of sign change for functions and distributions and a  history  of results. We will use Theorem~1.2 from \cite{MP_det_mes}.
To formulate this result, we need the following

\medskip\noindent
{\it Definition 4 }(\cite{MP_det_mes}). For a finite real measure $\mu$ on $\R$ we denote by $\mu_+$  and $\mu_-$ its mutually
singular positive and negative parts, $\mu=\mu_+-\mu_-,\mu_+\bot\mu_-.$ If $A$ and $B$ are
two disjoint closed subsets of $\R$, let $M_a(A, B)$ be the class of all finite real measures
$\mu$ with a spectral gap $[-a, a]$ such that supp$(\mu_+)\subset A$  and supp$(\mu_-)\subset B$. We
define the gap characteristic of a pair of closed subsets $A$ and $B$ of $\R$ as $$
G(A, B) := \sup\{a \geq 0 : M_a(A, B)\ne\emptyset\}.$$

\begin{theorem}[Theorem 1.2, \cite{MP_det_mes}]\label{tpol} For any closed sets $A, B \subset\R$,
$$G(A, B) = \frac{1}{2} \sup\{d: \exists  d-\mbox{{\rm uniform sequence} } \{\gamma_n\}=$$$$\{...<\gamma_{n-1}<\gamma_n<\gamma_{n+1}<...\}, \{\gamma_{2n}\} \subset  A,
\{\gamma_{2n+1} \}\subset  B\}.
$$\end{theorem}

See e.g. Definition 2.4 in \cite{MP_det_mes} for a definition of $d$-uniform sequence. For our application it is important that  every $d$-uniform sequence is  $d$-regular.

An immediate corollary of Theorem \ref{tpol} is that for every real function $u\in$\pw\! and  every $d<2\rho$ there is a real increasing $d$-regular sequence  $\Gamma=\{\gamma_k,k\in\Z\}$ such that $u(\gamma_{2k})\geq0$ and $ u(\gamma_{2k+1})\leq0, k\in\Z$.
Let Sign$(u)$ denote the set of real points where $u$ changes the sign. Clearly, between any two points $\gamma_n,\gamma_{n+1}$ there is a
point of  Sign$(u)$. It is also clear that Sign$(u)\subset Z_r(u)$. Using Definition 3 (i), we get

\begin{corollary}\label{c1}
Let $u\in PW_S^1$ be a real function. Then $D_{BM}(${\rm Sign}$(u))\geq 2\rho.$
\end{corollary}

\subsection{Proof of Proposition \ref{p1}}

\begin{lemma}\label{lemma}
Assume $\Gamma$ is a real sequence and $\Delta\subset\Gamma$ is a subsequence. Then
$$
D^{BM}(\Gamma\setminus\Delta)\leq D^{BM}(\Gamma)-D_{BM}(\Delta).
$$
\end{lemma}

\noindent {\bf Proof}.  If $D_{BM}(\Delta)=0$ then the statement is obvious.

Assume $D_{BM}(\Delta)>0$. Fix any  $\epsilon$ such that $0 < \epsilon < D_{BM}(\Delta)$. By Definition~3 (iii), every sequence of disjoint intervals $I_n$ such that
$$
\#(\Gamma\cap I_n)\geq (D^{BM}(\Gamma)+\epsilon)|I_n|,\quad \forall n,
$$is short.
Also, by Definitions 2 and 3 (i), every sequence of disjoint intervals $I_n$ such that
$$
\#(\Delta\cap I_n)\leq (D_{BM}(\Delta)-\epsilon)|I_n|,\quad \forall n,
$$is short.
Therefore, there is no long sequence of disjoint intervals $I_n$ such that for every $n$ one has
$$
\#((\Gamma\setminus\Delta)\cap I_n)=\#(\Gamma\cap I_n)-\#(\Delta\cap I_n)\geq (D^{BM}(\Gamma)-D_{BM}(\Delta)+2\epsilon)|I_n|.
$$
Now, Lemma \ref{lemma}  follows from Definition 3 (iii).

\medskip\noindent{\bf Proof of Proposition \ref{p1}}.
 Fix a function $f\in PW_S^1.$ Write $f(z)=u(z)+iv(z)$, where
$$
u(z):=(f(z)+\overline{f(\bar z)})/2,\quad v(z):=(f(z)-\overline{f(\bar z)})/(2i).
$$Clearly, $u$ and $v$ are real functions (on $\R$) which due to the symmetry of $S$  belong to $PW_S^1$.
At least one of these functions is nontrivial. We may assume that $u\ne0$. Clearly, $\La(u)=\La(f)$ and $\Omega(u)\supset\Omega(f)$. Hence, by Theorem~\ref{bm}, it suffices to check that $D^{BM}(\La(u)\cup\Omega(u))\leq 2(\sigma-\rho).$

Observe that Sign$(u)\cap \La(u)=\emptyset.$ Next, if a point $x$ is an element of both the set Sign$(u)$ and  the multiset $\Omega(u)$, since $u$ changes the sign at $x$ and since every point in $\Omega(u)$ has an even multiplicity $2n(x)$,
the multiplicity of zero at $x$ is exactly $2n(x)+1$. Therefore, $\La(u)\cup\Omega(u)\cup\,$Sign$(u)\subset Z(u)$. Using Lemma \ref{lemma} and Corollary \ref{c1}, we conclude that $$D^{BM}(\La(u)\cup\Omega(u))\leq D^{BM}(Z(u))- D_{BM}(\mbox{{\rm Sign}}(u))\leq 2\sigma- 2\rho,$$which proves the proposition.

\subsection{Proof of necessity part of Theorem \ref{t1}}

The necessity part of Theorem \ref{t1} follows from
\begin{proposition}\label{p2} {\rm (i)} Assume $f\in\,$\ext. Then $f$ satisfies {\rm (\ref{c0})}, {\rm (\ref{sp_cond}),}  and {\rm (\ref{la_cond})}.

{\rm (ii)} Assume $f\in\,$\expo. Then $f\in\,$\ext\! and $f$ satisfies {\rm (\ref{exp2})}, and {\rm (\ref{exp1})}.
\end{proposition}

\noindent
{\bf Proof}. The proof is almost identical  to the proof of the necessity part of Theorem 3 in \cite{D_PW}.

For the reader's convenience, we write the prove of (i).   It suffices to check that $f$ satisfies  (\ref{sp_cond})  and (\ref{la_cond}).

Firstly, if (\ref{sp_cond}) is not fulfilled, one may put $h(x) = \sin(\varepsilon x)$. Then Sp$(hf)\subset\,$Sp$(f)+\{-\varepsilon,\varepsilon\}\subset S$,  for any sufficiently small $\varepsilon > 0$, and so  $h f \in PW_S^1$. Then by  Lemma~\ref{h_lemma}, $f\notin$\ext.

Secondly, if condition (\ref{la_cond}) is not satisfied,  then there is $\lambda \in \Cc\setminus\R$ such that
$
f_+(\lambda) = f_-(\lambda) = f_+(\bar \lambda) = f_-(\bar \lambda) = 0,
$
where $f = f_+ + f_-$, Sp$(f_-) \subset [-\sigma, -\rho]$, and Sp$(f_+) \subset [\rho, \sigma]$.  Consider the real bounded function
$$
h(x) := \frac{1}{ x - \bar \lambda} + \frac{1}{x - \lambda}.
$$
One may easily check that the functions $f_- h, f_+ h$ are entire and that Sp$(f_- h) \subset [-\sigma, -\rho]$ and Sp$(f_+  h) \subset [\rho, \sigma]$. This shows that $f h \in PW^1_{S}$, and by Lemma~\ref{h_lemma}, $f \notin$ Ext$(PW^1_S)$.

\subsection{Proof of sufficiency of part (A)}

 Assume that $f$ satisfies the assumptions in part (A) of Theorem \ref{t1}. Assume also that there is a real bounded function $h$ such that  $g:=fh\in PW_S^1$. To prove the sufficiency, by Lemma \ref{h_lemma} it suffices to show that $h= const.$

We have $h=g/f$. There are two possibilities:
\begin{enumerate}
    \item[(i)] $h$ is an entire function;
    \item[(ii)] $h$ has poles.
\end{enumerate}

\medskip
(i) Assume that $h$ is an entire function.  Recall that $h$ is  real and bounded on the real line. Therefore,
$d:=$Type$(h)=$Type$_-(h)=$Type$_+(h)$.
Also, by Lemma~\ref{ko} (iii),  we see that $d= $Type$(g)-$Type$(f)\leq 2\pi(\sigma-\rho)<2\pi\rho$.

 If $d=0$, then Lemma \ref{l001} (iii)  implies $h=const$.

Assume $d>0.$ If $\{-\sigma\}\in$Sp$(f)$, then Type$_+(f)=2\pi\sigma$. Using Lemma \ref{ko} (ii), we get \begin{equation}\label{koi0}d=\mbox{{\rm Type}}_+(h)=\mbox{{\rm Type}}_+(g)-\mbox{{\rm Type}}_+(f)\leq2\pi(\sigma-\sigma)\leq0,\end{equation} which means that $h$ is of zero type. Contradiction. Similarly, condition $\{\sigma\}\in$Sp$(f)$ leads to a contradiction.

Assume that $\{-\rho\}\in$Sp$(f)$. Write $g=g_-+g_+$, where Sp$(g_-)\subset[-\sigma,-\rho]$ and Sp$(g_+)\subset[\rho,\sigma]$. One may check that $g_-=hf_-$ and $g_+=hf_+$.
However, using Lemma \ref{ko} (ii), we get \begin{equation}\label{koi}-2\pi\rho\geq\mbox{{\rm Type}}_-(g_-)=\mbox{{\rm Type}}_-(f_-)+\mbox{{\rm Type}}_-(h)=2\pi(-\rho+d).\end{equation} Contradiction.

Alternatively, one may use the Titchmarsh convolution theorem for distributions
with compact support and the fact that $d=$Type$(h)<2\pi\rho$ to see that  $\sup\,$Sp$(f_-h)=-\rho+d/(2\pi)<0$. Similarly, $\inf\,$Sp$(f_+h)>0$. Therefore, the spectrum of $g=fh$ cannot lie on $S$. Contradiction.


\medskip

(ii)  Now, assume that $h=g/f$ is a meromorphic function which has poles. By Lemma \ref{l001} (i), the poles of $h$ may only lie on the set $\La(f).$

Assume  $\{-\sigma\}\in$Sp$(f)$.
Choose any point $\la$ which is a pole of  $h$. Since $h$ is real, then  $\bar\lambda$ is also a pole of $h$. Let us assume that ${\rm Im} \, \lambda < 0$. Let $m\geq 1$ denote the order of pole $\lambda$.

By Corollary~\ref{c00},  $\rho':=R(\La(f))<\rho$.  This means that for every $r$ satisfying $\rho'<r<\rho$ there is a  function $\varphi  \in PW^1_{[-r,r]}$ such that the function $h\varphi$ is entire and does not vanish at $\la$. Hence, $\la$ is a zero point of $\varphi$ of multiplicity $m$. We may assume that $\varphi$ is {\it real} on $\R$ (Otherwise, we consider $\varphi(z)+\overline{\varphi(\bar z)}$ or $i(\varphi(z)-\overline{\varphi(\bar z)})$). Since $\{-\sigma\}\in$Sp$(f)$ and Sp$(g)\subset[-\sigma,\sigma]$, we see that the function $h\varphi=g\varphi/f$ is at most of type $2\pi r$ in the upper half-plane. Since this function is real, the same is true for the lower half-plane. We conclude that $h\varphi\in PW_{[-r,r]}^1$.

 Set $\psi(z):=\varphi(z)/(z-\la).$
This function also belongs to $PW_{[-r,r]}^1.$ We have
$fh\psi=g\psi$. By the Titchmarsh convolution theorem, both functions $g\varphi$ and $g\psi$ have spectral gap $(-\rho+r,\rho-r).$ On the other hand, the Fourier transform of the function $fh\psi=fh\varphi/(z-\la)$ is equal to
$F\ast Q\ast e_\la$, where $F$ and $Q$ denote the Fourier transform of $f$ and $h\varphi$, respectively, and
\begin{equation}\label{dl}
e_\la(t)=\left\{\begin{array}{ll} -2 \pi i e^{- 2\pi i \la t},& t\ge0\\ 0& t<0
\end{array}\right.
\end{equation}is the Fourier transform of $1/(z-\la)$. Write $f=f_-+f_+$ and denote by  $F_-$ the Fourier transform of $f_-$ which vanishes outside  $[-\sigma,-\rho]$. Since $e_\la(t)=0$  for $t<0$ and $Q$ vanishes outside $(-r,r)$, one may check that for $|t|<\rho-r$ we have
$$F\ast Q\ast e_\la(t)= F_-\ast Q\ast e_\la(t)=-2\pi i f_-(\la)(h\varphi)(\la)e^{-2\pi i \la t}.$$ Recall, by (\ref{la_cond}), that $f_-(\la)\ne0$, and from above we know that
$(h\varphi)(\la)\ne0$. We conclude that $F\ast Q\ast e_\la(t)\ne0$ on $(-\rho+r,\rho-r)$.
Contradiction. Similarly, condition $\{\sigma\}\in$Sp$(f)$ leads to a contradiction.

Finally, we may assume that  $\{-\rho\}\in$Sp$(f)$, while none of the points $\{\pm \sigma\}$ lies in Sp$(f)$. Then
\begin{equation}\label{r}
r':=\max\{r>0:\{ \pm r\} \cap \mbox{Sp}(f)\} \neq \emptyset\}\in (\rho,\sigma).
\end{equation}
We see that $f\in PW^1_{[-r',-\rho]\cup[\rho,r']}.$ By Proposition \ref{p1},  $R(\La(f))\leq r'-\rho$. This means that for every $r>r'-\rho$ there is a real function $\varphi\in PW_{[-r,r]}^1$ such that $h\varphi=g\varphi/f$ is an entire function. Since $\rho>\sigma/2$, we may assume that $r<r'+\rho-\sigma $. Then it is easy to check that the type of $h\varphi$ is at most $2 \pi (\sigma+r-r')$, a number which is less than $ 2 \pi \rho$.
The rest of the proof is as in the previous step, which leads to a contradiction.

\subsection{Proof of sufficiency of part (B)}

Assume that $f\in$\ext, $f$ satisfies (\ref{exp2}), (\ref{exp1}) and  $fh=g\in$\pw\! for some nonnegative  function $h$. To prove the theorem, by Lemma \ref{h_lemma}, we have to show that $h=const.$

The proof is by contradiction. We assume that $h\ne const.$ Then $h=g/f$ is either (i) nonconstant entire or (ii) meromorphic function.

\medskip

(i) Assume $h$ is a nonconstant entire function.   The proof is pretty similar to the proof of part (i) in Sec. 2.6.

By Lemma \ref{ko}, $f\in$Cart and is of exponential type $d:= $Type$(g)-$Type$(f)<2\pi\rho$.

Assume  $d=0$. This means that $h$ is an entire function of zero type. By (\ref{exp2}), the function $g=fh$ cannot lie in \pw\!. Contradiction.

Assume $d>0$. If $\{-\sigma\}\in$Sp$(f)$, one may  use (\ref{koi0}) to arrive to contradiction. Similarly, condition $\{\sigma\}\in$Sp$(f)$ leads to a contradiction.

 If $\{-\rho\}\in$Sp$(f)$, then (\ref{koi}) leads to a contradiction. Similarly, condition $\{\rho\}\in$Sp$(f)$ leads to a contradiction.

\medskip
(ii) Assume that $h$ has poles.
Denote by $P(h)$ the set of poles of $h$ counting multiplicities.  Every pole of $h$ is a zero of $f$. Since $h$ is nonnegative on $\R$, every real pole $\widetilde{x}$ of $h$ is of an even order $2n(\widetilde{x})$. Also, if $\lambda\in\Cc\setminus\R$ is a purely complex pole, then $\bar\lambda$ is also a pole of $h$. Therefore, $P(h)\subset\La(f)\cup\Omega(f)$.

If $h$ has a purely complex pole $\lambda\in\La(f)$, we  proceed exactly as  in the part (ii) of the proof in Sec. 2.6: We use  Corollary \ref{c0} to find two functions $\varphi(z)$ and $\psi(z)=\varphi(z)/(z-\lambda)\in PW^1_{[-r,r]}$, $0<r<\rho,$ such that both functions $fh\varphi$ and  $fh\psi$ belong to $ PW^1_{[-\rho+r,\rho-r]}$. This means that $fh\psi$ has  a spectral gap at the origin. However, as in Sec. 2.6 we show that this is not possible. Contradiction.

Assume that $h$ has only real poles. Choose any point $x_0\in P(h)$. 
By (\ref{omega}),  one of the conditions hods: either $f_-(x_0)\ne0$ or  $f_-(x_0)=0, f_-'(x_0)\ne0.$

Let us assume that $f_-(x_0)=0, f_-'(x_0)\ne0.$
 By Corollary \ref{c0}, there is number $r, 0 < r <\rho,$ and a function $\varphi\in PW_{[-r,r]}^1$ such that the function
$h\varphi$ is entire and does not vanish on $P(h)$. 
Set $\psi (z) := \varphi(z)/(z-x_0)^2\in PW^1_{[-r,r]}$. 

\begin{claim} The Fourier transform of
the function $fh\psi=g\psi  = g\varphi/(z-x_0)^2$ is equal to $ Q \ast \widetilde{e}_{x_0}$, where   $Q$ is
the Fourier transform of $g\varphi$ and $$\widetilde{e}_{x_0}(t): =\left\{ \begin{array}{ll}
-(2\pi)^2 te^{-2\pi i x_0t} & t\geq0\\ 0& t<0.
\end{array}\right.
$$\end{claim}

Indeed, by (\ref{dl}),  the Fourier transform of $1/(z-\lambda)^2,$ Im$\lambda<0,$ is equal to $-(2\pi)^2te^{2\pi i \lambda t}$ for $t>0$, and it vanishes for $t<0$. 
To prove the claim, one may consider $\lambda=x_0-i\varepsilon, \varepsilon>0$, and let $\epsilon\to0$.
We leave the details to the reader.

Since $g\varphi\in PW^1_{[-\sigma - r,-\rho+r]\cup[\rho-r,\sigma+r]}$, we have $Q(t)=0, |t|\leq \rho-r$.
Write $g\varphi=(g\varphi)_-+(g\varphi)_+$ and $Q=Q_-+Q_+$, where  $$Q_\pm:=\widehat{(g\varphi)_\pm}=\widehat{(fh\varphi)_\pm}=\widehat{f_\pm h\varphi}.$$Clearly, $Q_-$ vanishes for $t\geq -\rho+r$ and  $Q_+$ vanishes for $t\leq \rho-r$.
As in the proof in Sec. 2.6,  for $t\in (-\rho+r,\rho-r)$ we get
$$
 Q\ast \widetilde{e}_{x_0}(t) =  Q_-\ast \widetilde{e}_{x_0}(t) =-2\pi i( f_-h\varphi)'(x_0)e^{-2\pi i x_0t}.
$$Since $( f_-h\varphi)'(x_0)=f_-'(x_0) (h\varphi)(x_0)\ne0$, we conclude that $ Q \ast \widetilde{e}_{x_0}(t)\ne0$  on $(-\rho+r,\rho-r)$. Contradiction.

If condition $f_-(x_0)\ne0$ holds, we set $\psi (z) := \varphi(z)/(z-x_0)$ and proceed as above.

\subsection{Extensions of Theorem \ref{t1}}
An inspection of the  proof in Sec. 2.6  shows that it  admits the following two extensions, which we state without proof.

Firstly,  the  sufficiency part of Theorem \ref{t1} (A) admits extension  to the spectra $S$ with short gap:

\begin{theorem}[Short gap]\label{t7} Assume $0<\rho<\sigma/2$ and that  a function  $f\in$\pw satisfies {\rm(\ref{c0})}, {\rm(\ref{sp_cond})},  {\rm(\ref{la_cond})},  and \begin{equation}\label{0}R(\La(f))<\rho+r'-\sigma,\end{equation} where $r'$ is defined in {\rm(\ref{r})}.
Then $f\in${\rm \ext}.
\end{theorem}

The restriction (\ref{0}) makes sense only if $\rho+r'>\sigma$. See the example in Sec.~5 of a function $f\notin$\ext\! satisfying {\rm(\ref{c0})}, {\rm(\ref{sp_cond})},  {\rm(\ref{la_cond})},  and $\La(f)=\emptyset.$ Note that if at least one of the end-points $\{\pm\sigma\}$ lies in Sp$(f)$,
then (\ref{0}) becomes $R(\La(f))<\rho$. From Theorems~\ref{t2} and \ref{t7} one gets

\begin{corollary}[Short gap]\label{carb} Assume that $0<\rho<\sigma/2$.

{\rm(i)} If $f\in$\pw satisfies {\rm(\ref{c0})}, {\rm (\ref{d1})},  {\rm(\ref{la_cond})} and $R(\La(f))<\rho$, then $f\in${\rm \ext}.

{\rm (ii)} There exists  $f\in$\pw satisfying  {\rm(\ref{c0})}, {\rm (\ref
{d1})},   {\rm(\ref{la_cond})}
and  $R(\La(f))>\rho$ such that  $f\notin${\rm \ext}.
\end{corollary}

Secondly, Theorem \ref{t1} remains valid for the union of two arbitrary intervals, provided the gap between the intervals is large:
\begin{equation}\label{sets}
S=[a,b]\cup[c,d],\quad a<b<c<d,\ c-b>(d-a)/2.
\end{equation}
Similarly to above, write  $f=f_-+f_+,$  where Sp$(f_-)\subset[a,b]$, Sp$(f_+)\subset[c,d]$.

\begin{theorem}\label{tlast} Let $S$ satisfy {\rm (\ref{sets})}. Then the set {\rm \ext} consists of the functions $f\in$\pw\! satisfying {\rm (\ref{c0})}, {\rm (\ref{la_cond})} and $$
\mbox{at least one of the points  }\{a\}, \{b\}, \{c\}, \{d\} \mbox{ belongs to {\rm Sp}}(f),
$$

\end{theorem}

Similar extensions of part (B) of Theorem \ref{t1} are also true.

\section{Short gap: Proof of Theorem \ref{t2}}
 Assume $\sigma>2\rho$ (short gap).
By a change of variables, we may assume that $\sigma=2+\epsilon$ and $\rho=1-\epsilon,$ i.e. we consider the spectra
$$
S_\epsilon:=[-2-\epsilon,2+\epsilon]\setminus(-1+\epsilon, 1-\epsilon),\quad \mbox{where } 0<\epsilon<1.
$$


\subsection{Non-extreme function with $\La(f) = (\Z+i)\cup(\Z-i)$}

 Set
 \begin{equation}\label{phi}
 \varphi_\epsilon(z):=i\left(\frac{\sin(\pi \epsilon z)}{\pi \epsilon z}\right)^2
 \end{equation} and
\begin{equation}\label{f}
f(z):=4(\cos(2\pi z)-\cos(2\pi i))\left(\cos(2\pi z)+\frac{1}{2\cos(2\pi i)}\right) \varphi_\epsilon(z).
\end{equation}

\begin{lemma}\label{l1}
We have $f\in PW_{S_\epsilon}^1$, $f$ satisfies {\rm(\ref{la_cond})}, $\{-2-\epsilon\}\in${\rm S}p$(f)$, and  $\La(f)=\Z \pm i$.
\end{lemma}

We skip the simple proof.

By Lemma \ref{l1}, part (i) of  Theorem \ref{t2} is an immediate consequence of the following

\begin{lemma}\label{lt}
 $f\!\notin\,${\rm Ext}$(S_\epsilon)$.
\end{lemma}

To prove this  we will need a number of  lemmas.

\begin{lemma}\label{l3}
There are coefficients $c_n\in\Cc$ such that $|c_n|\leq (1+|n|)^{-4}, n\in\Z,$ and
$$
\sum_{n\in \Z}\left(c_ne^{2\pi i (n+i)t}+\bar c_ne^{2\pi i (n-i)t}\right)=0,\quad |t|\leq 1-\epsilon.
$$
\end{lemma}

\noindent {\bf Proof}. We start with

\medskip\noindent
\begin{claim}\label{cl00} Every function $Q\in L^2(-1,1)$ can be written in the form $Q(t)=F_1(t)+F_2(t),$ where
$$
F_1(t):=\sum_{n\in\Z}a_ne^{2\pi i (n+i)t}=e^{-2\pi t}\sum_{n\in \Z}a_ne^{2\pi i nt},\quad \{a_n\}\in l^2(\Z),$$$$ F_2(t):=\sum_{n\in\Z}b_ne^{2\pi i (n-i)t}=e^{2\pi t}\sum_{n\in \Z}b_ne^{2\pi i nt},\quad \{b_n\}\in l^2(\Z).
$$\end{claim}

Clearly, we have
$$
F_1(t+1)=e^{-2\pi}F_1(t),\quad F_2(t+1)=e^{2\pi}F_2(t), \quad t\in\R.
$$
Therefore, to prove Claim \ref{cl00} we have to find $F_1,F_2$ so that
$$
F_1(t)+F_2(t)=Q(t),\  e^{-2\pi}F_1(t)+e^{2\pi}F_2(t)=Q(t+1),\quad -1<t<0.
$$This is equivalent to
$$
F_1(t)=\frac{e^{2\pi}Q(t)-Q(t+1)}{e^{2\pi}-e^{-2\pi}}:=Q_1(t),\ F_2(t)=\frac{Q(t+1)-e^{-2\pi}Q(t)}{e^{2\pi}-e^{-2\pi}}:=Q_2(t),$$ where $-1<t<0.
$
Clearly, $Q_j\in L^2(-1,0), j=1,2.$ Finally, since the trigonometric system $\{e^{2\pi i nt},n\in\Z\}$
forms an orthonormal basis for $L^2(-1,0),$ it is clear that 
we may choose coefficients $a_n,b_n$ so that for $0<t<1$ we have
$$
e^{2\pi t}F_1(t)=\sum_{n\in\Z}a_ne^{2\pi i nt}=e^{2\pi t}Q_1(t),\quad e^{-2\pi t}F_2(t)=\sum_{n\in\Z}b_ne^{2\pi i nt}=e^{-2\pi t}Q_2(t),
$$which proves the claim.

Choose any continuous function $Q$ on $[-1,1]$ which vanishes for  $|t|\notin (1-\epsilon/2,1)$. By the claim above, there are  $l^2$-sequences $a_n,b_n$ such that
$$
Q(t)=\sum_{n\in\Z} \left(a_n e^{2\pi i (n+i)t}+b_ne^{2\pi i (n-i)t}\right)=0,\quad \mbox{a.e. for } |t|<1-\epsilon/2.
$$
Take any non-trivial smooth function $\psi\in C^\infty(\R)$ which vanishes outside $(-\epsilon/2,\epsilon/2)$. Then its Fourier transform $|\hat\psi(x\pm i)|$ tends to zero as $|x|\to\infty$ faster than  $|x|^{-N}$, for any $N>0$. Consider the convolution $\psi\ast Q$. Clearly,
$$
\psi\ast Q(t)=\sum_{n\in\Z} \left(\hat\psi(n+i)a_n e^{2\pi i (n+i)t}+\hat\psi(n-i)b_ne^{2\pi i (n-i)t}\right)=0,\quad |t|<1-\epsilon.
$$
Set $c_n:=\delta(\hat\psi(n+i)a_n+\overline{\hat\psi(n-i)b_n}),$ where $\delta>0$ is a small number. Then we have $|c_n|\leq (1+|n|)^{-4}$, $n\in\Z,$
and $$
 \delta(\psi\ast Q(t)+\overline{\psi\ast Q(-t)})=\sum_n \left(c_n e^{2\pi i (n+i)t}+\bar c_ne^{2\pi i (n-i)t}\right)=0,\quad |t|<1-\epsilon.
$$


\begin{lemma}\label{l4}
Denote by $\Phi$ the Fourier transform of the function $\varphi_{\varepsilon}$ in {\rm (\ref{phi})}. Then the Fourier transform $F(t)$ of $f$ in {\rm (\ref{f})} is given by
$$
F(t)=\Phi(t-2)+A\Phi(t-1)+A\Phi(t+1)+\Phi(t+2),\quad t\in\R,
$$where $A:=-2\cos(2\pi i)+1/\cos(2\pi i)$.
\end{lemma}

We skip the simple proof.

\begin{lemma}\label{l5}
The Fourier transform $F_n(t)$ of the function $$f(z)\left(\frac{c_n}{z-n-i}+\frac{\bar c_n}{z-n+i}\right)$$ satisfies $F_n(t)=0, |t|\geq 2+\epsilon$, and
$$
F_n(t)=K\left(c_n\varphi_\epsilon(n+i)e^{-2\pi i (n+i)t}+ \overline{c_n\varphi_\epsilon(n+i)}e^{-2\pi i (n-i)t}\right),\quad |t|\leq 1-\epsilon,
$$where
$$
K:=2\pi i\left(Ae^{-2\pi}+e^{-4\pi}\right).
$$
\end{lemma}

\noindent {\bf Proof}.
Clearly,
$$
\frac{1}{z-n-i}=\mathfrak{F}^{-1} u_n(z),\quad \frac{1}{z-n+i}=-\mathfrak{F}^{-1} {v_n}(z),
$$where $\mathfrak{F}^{-1}$ denotes the inverse Fourier transform,
$$
u_n(t):=\left\{\begin{array}{ll} 0& t>0\\ 2\pi i e^{-2\pi i (n+i)t} & t<0\end{array}\right.
$$and
$$
v_n(t):=\left\{\begin{array}{ll} 2\pi i e^{-2\pi i (n-i)t} & t>0\\ 0& t<0\end{array}\right.
$$
Therefore
$$
\varphi_\epsilon(z)\frac{1}{z-n-i}=\mathfrak{F}^{-1} \left(\Phi\ast u_n\right)(z).
$$

Since Sp$(\varphi_\epsilon)\subset[-\epsilon,\epsilon]$ and $\varphi_\epsilon(-z)=\varphi_\epsilon(z)$, one may check that
$$
\Phi\ast u_n(t)=\left\{\begin{array}{ll} 0& t>\epsilon\\ 2\pi i \varphi_\epsilon(n+i)e^{-2\pi i (n+i)t}& t<-\epsilon\end{array}\right.
$$
and
$$
\Phi\ast v_n(t)=\left\{\begin{array}{ll}  2\pi i \varphi_\epsilon(n-i)e^{-2\pi i (n-i)t}& t>\epsilon\\0& t<-\epsilon\end{array}\right.
$$

From  above and Lemma \ref{l4}, for
$|t|\leq 1-\epsilon$ we obtain
$$
F_n(t)=K\left(c_n\varphi_\epsilon(n+i)e^{-2\pi i (n+i)t}-\bar c_n\varphi_\epsilon(n-i)e^{-2\pi i (n-i)t}\right).
$$
By (\ref{phi}), $\varphi_\epsilon(n-i)=-\overline{\varphi_\epsilon(n+i)}$.

Finally, since $f(n-i)=0,$ we get for $t\geq 2+\epsilon$
$$
F_n(t)=F\ast(c_nu_n-\bar c_n v_n)(t) =-\bar c_n F\ast v_n(t)=$$$$-\bar c_n f(n-i)e^{-2\pi i (n-i)t}=0.
$$ Similarly, $F_n(t)=0$ for $t\leq-2-\epsilon,$ which completes the proof.

\medskip


\noindent {\bf Proof of Lemma \ref{lt}.}
Let $c_n$ be the coefficients in Lemma \ref{l3}. Set
$$
h(z):=\sum_{n\in\Z}\left(\frac{c_n}{ \varphi_\epsilon(n+i)}\frac{1}{z-n-i}+\frac{\bar c_n}{ \overline{\varphi_\epsilon(n+i})}\frac{1}{z-n+i}\right).
$$Then $h$ is real on $\R$. By the estimate on $c_n$ in Lemma \ref{l3}, $h$ is bounded on $\R$.

Finally, it follows from Lemmas \ref{l3} and \ref{l5} that the Fourier transform of  $fh$
vanishes for $t\in [-1+\epsilon,1-\epsilon]$ and for $|t|\geq 2+\epsilon.$ Therefore, $fh\in PW_{S_\epsilon}^1$.
 Lemma~\ref{h_lemma} yields  $f\notin$Ext$(PW_{S_\epsilon}^1)$.

 \subsection{Non-exposed function with $\Omega(f) = \Z$}
 To prove part (ii) of Theorem \ref{t2}, for every $\varepsilon, 0<\varepsilon<1,$ we have to find a function $f\in$Ext$(PW_{S_\varepsilon^1})$ satisfying (\ref{exp2}), (\ref{exp1}) and
 such that $f\notin$Exp$(PW_{S_\varepsilon}^1)$.

 Set
 $$
\varphi(z):=\frac{\sin (2\pi\varepsilon z)}{(\varepsilon z)^2-1}, \  f(z):=4c\sin^2(\pi z)\sin(2\pi z)\varphi(z),
$$where $c$ is a constant such that $\|f\|_1=1$. Clearly, $\varphi\in PW_{[-\varepsilon,\varepsilon]}^1$ and $\{\pm\varepsilon\}\in$Sp$(\varphi)$.

By Lemma~\ref{h_lemma}, part (ii) of Theorem \ref{t2} follows immediately from

\begin{lemma}\label{cl} $f$ satisfies {\rm (\ref{exp2}), (\ref{exp1})}, $f\in${\rm Ext}$(PW_{S_\varepsilon}^1)$, and  $f(z)/\sin^2{\pi z}\in PW_{S_\varepsilon}^1.$
\end{lemma}

\noindent{\bf Proof}. Write
\begin{equation}\label{sincos}4\sin^2(\pi z)\sin(2\pi z)= 2\sin(2\pi z)-\sin(4\pi z).
\end{equation}Therefore,  Sp$(\sin^2(\pi z)\sin(2\pi z))=\{-2,-1,1,2\}$. It  follows that $f\in PW_{S_\varepsilon}^1.$ It is also clear that
$\{-2-\epsilon\}\in$Sp$(f)$ and that $f$ has only real zeros. Therefore, by part (A) of Theorem \ref{t0} and inclusion~(\ref{ex_inclusion}),  $f\in$Ext$(PW_{S_\varepsilon}^1)$.

Next, from (\ref{sincos}) we get
 $f=f_-+f_+,$ where
$$
f_-(z)= -c\left(e^{-4\pi iz}+2e^{-2\pi i z}\right)\varphi(z), \quad f_+(z)= c\left(e^{4\pi iz}+2e^{2\pi i z}\right)\varphi(z).
$$It follows that $f$ satisfies (\ref{exp1}).

Now, assume $w(x)$ is an entire nonnegative function of zero type such that $wf\in L^1(\R).$ To prove (\ref{exp2}), we have to show that $w=const$.  Since $wf$ is an entire function of exponential type $4\pi +2\epsilon,$ it satisfies (\ref{in}) with some constant $C$ and $d=4+2\varepsilon.$ One can easily check that this implies
$$
|w(z)|\leq C_1(1+|z|^2),\quad z=x+iy, x\in\R, |y|\geq 1, C_1=const.
$$
Since $w|_\R\geq0,$ it easily follows that either $w=const$ or $w(z)=(z-a)^2$ for some $a\in\R$. The latter is not possible since it implies $wf\notin L^1(\R)$. We conclude that $w=const.$

Finally, since $$\mbox{{\rm Sp}}(\sin(2\pi z)\varphi(z))=[-1-\varepsilon,-1+\varepsilon]\cup[1-\varepsilon,1+\varepsilon]\subset S_\varepsilon,$$ we see that
$$
\frac{f(z)}{\sin^2{\pi z}}=4c\sin(2\pi z)\varphi(z)\in PW_{S_\varepsilon}^1.
$$

\section{Short Gap: Proof of Theorem \ref{t3}}

\subsection{Extreme function with $\La(f)=(\Z+i)\cup(\Z-i)$}
Here we  prove part (i) of Theorem \ref{t3}.

 By a change of variables, we may assume that $\rho=1-\de,\sigma=4+\de$, where $0<\de<1$, i.e.
\begin{equation}\label{s}
S=[-4-\de,4+\de]\setminus(-1+\de,1-\de),\quad 0<\de<1.
\end{equation}

It follows from Sec. 3.1 that there is a function $f\in$\pw satisfying {\rm(\ref{c0})}, {\rm(\ref{sp_cond})}, {\rm(\ref{la_cond})}, $\La(f)=\Z\pm i$ and such that $f\notin$\ext. To prove part (i) of Theorem \ref{t3}, we  will construct a function $f\in$\pw satisfying (\ref{c0}), (\ref{sp_cond}), (\ref{la_cond}), $\La(f)=\Z\pm i$ and such that $f\in$\ext.

For brevity, in this section, we set $A:= \cos(2 \pi i)$ and $c = \frac{3A - 4 A^3}{2A^2 - 1}$.
We start with the following
\begin{lemma}\label{tau}
The function $\tau$ defined by
\begin{equation}\label{tau_def}
    \tau(z): =2\left( \cos(6 \pi z) + c \cos(4 \pi z)\right)
\end{equation}
has the following properties
\begin{equation}\label{tau_la}
\Lambda(\tau) = Z_c(\tau)=\Z \pm i, 
\end{equation}
and for every $\varphi \in L^1(\R)$ we have
\begin{equation}\label{phi_1}
    \widehat{\tau \varphi}(t ) = \Phi(t-3) + c \Phi(t-2) + c\Phi(t+2) + \Phi(t+3),
\end{equation}
where $\Phi:=\hat{\varphi}$.
\end{lemma}

{ \bf Proof.}
One may check that $\tau$ can be rewritten as
$$
\tau(z) =  2( \cos(2 \pi z) - A)\left( 2 \cos(4 \pi z) + \frac{2A}{2A^2-1} \cos(2 \pi z) + \frac{1}{2A^2-1} \right).
$$

Using Rouché's theorem (see details in \cite{ou}, sec. 5), one may check that the  second factor above has only real zeros. The first factor vanishes exactly on $\Z\pm i$. Hence,  $\Lambda(\tau) = \Z \pm i$.
Finally, (\ref{phi_1}) easily follows from the definition of  $\tau$ in (\ref{tau_def}).


\begin{remark} Write $\tau=\tau_-+\tau_+$. By direct computation, one can show that the functions
\begin{equation}\label{tau12_def}
   \tau_+(z) = e^{6 \pi i z} + c e^{4 \pi i z} \quad \text{ and } \quad \tau_-(z) = e^{-6 \pi i z} + c e^{-4 \pi i z}
\end{equation}
do not have zeros  on $\Z \pm i$.
\end{remark}

Next, set $\varepsilon_n = \eta 4^{-n^2}$ and $\delta_n =\eta 3^{-n^2}, \, n \in \N,$ where $\eta>0$ is a small number. Set
\begin{equation}\label{pp}
    \psi_{\varepsilon}(z): = \prod \limits_{n = 1}^{\infty}\left(1 - \frac{z^2}{(n+\varepsilon_n)^2} \right) \quad \text{and} \quad
    \psi_{\delta}(z): = \prod \limits_{n = 1}^{\infty}\left(1 - \frac{z^2}{(n+\delta_n)^2} \right).
\end{equation}
Both functions  $ \psi_{\varepsilon}(z)$ and $ \psi_{\delta}(z)$ can be viewed as small `perturbations' of the function $$\frac{\sin(\pi z)}{\pi z}=\prod_{n=1}^\infty\left(1 - \frac{z^2}{n^2} \right).$$

Set
\begin{equation}\label{funfi}
\varphi(z) := \psi_{\delta}(z+i) \psi_{\varepsilon}(z-i) \frac{\sin(2\pi d z)}{z}
\end{equation}

\begin{lemma}For all small enough $\eta$ we have $\varphi\in PW^1_{[-\de,1+\de]}$ and
\begin{equation}\label{small_res_est}
    |\varphi(n \pm i)| \le  C2^{-n^2},\quad n\in\Z, \quad C=const.
\end{equation}\end{lemma}

To prove the estimate above, one can use  a standard trick in complex analysis: One compares the growth of $\psi_{\varepsilon}(z)$ and $\psi_{\delta}(z)$ with the growth of
$\sin(\pi z)$, see i.e. \cite{ou}, the proof of Lemma 4.29.
In particular, one can check that
 $\psi_{\varepsilon}\psi_{\delta}$ belong to $PW^1_{[-1,1]}$, provided $\eta$ is small enough. Therefore,
  Sp$(\varphi)\subset[-1,1]+[-\de,\de]=[-1-\de,1+\de]$, and so $\varphi\in PW^1_{[-1-\de,1+\de]}$.

 Set
\begin{equation}\label{funf}
    f(z) := \tau(z) \varphi(z).
\end{equation}
Then
$$
{\rm Sp}(f)\subset {\rm Sp}(\tau)+{\rm Sp}(\varphi)=\{-3,-2,2,3\}+[-1-\de,1+\de]=[-4-\de,4+\de]\setminus(-1+\de,1-\de).
$$By (\ref{s}), this implies $f \in$\pw. Moreover, from (\ref{tau_la}) and the definitions of the functions $f, \psi_{\delta}$, and $\psi_{\varepsilon}$ we see that $$ \Lambda(f) = \{\Z \pm i\}.$$

Finally, it remains to prove
\begin{lemma}\label{sist}
 $f\in$\ext.
\end{lemma}

\noindent
{\bf Proof}.
We will use the extreme point criterion formulated in Lemma~\ref{h_lemma}.
Assume that there is a real function $h$ from $L^{\infty}(\R)$ such that $fh = g$ and $g\in PW^1_S$. Our aim is to show that $h = const$. We will argue by contradiction: assuming that $h \neq const$, we will show that $g$ cannot have a spectral gap at the origin. This leads to a contradiction, since every function $g\in$\pw has spectral gap $[-1+\de,1-\de]$.

Similarly to Step (ii)  of the proof of Theorem \ref{t1}, we see that
$$
Z(f) \setminus \{\Z \pm i\} \subset Z(g).
$$
Hence, $h$ is a meromorphic function which may have simple poles only at some points from ${\Z \pm i}$.  By construction, $f$ is of exponential type $2 \pi(4+\de)$ and $g$ is of exponential type $\leq 2 \pi( 4+\de)$. If $h$ does not have poles, it is an entire function of zero type. Since it is also bounded on $\R$, it follows that $h\equiv const,$ which is a contradiction.
 Therefore, $h$  has poles.

Observe that the function $\cos(2\pi z)-\cos(2\pi i)$ is of exponential type $2\pi$ and it vanishes on $\Z\pm i$. Therefore, by (\ref{tau_la}) and (\ref{funf}),  the function
$$
q(z):=h(z)(\cos(2\pi z)-\cos(2\pi i))=\frac{g(z)(\cos(2\pi z)-\cos(2\pi i))}{f(z)}
$$is an entire function. It is real and bounded on $\R$.  Similarly to Sec. 2.5 above, we conclude   $q$ is of finite exponential type  $\leq 2\pi$.
In particular, it is bounded on every line in $\Cc$ parallel to $\R$.

Since $\varphi$ does not vanish on $\Z\pm i,$  $\varphi(z) h(z)$ has simple poles at each pole of $h$. Using  (\ref{small_res_est}) and (\ref{funf}) it is easy to check that the residues at these points
$$
    c_n : = \Res_{z = n + i} \varphi (z) h(z)=\Res_{z = n + i}\frac{ \varphi (z) q(z)}{\cos(2\pi z)-\cos(2\pi i)}=-\frac{\varphi(n+i)q(n+i)}{2\pi \sin(2\pi i)},$$$$ d_n : = \Res_{z = n - i} \varphi (z) h(z)=\frac{\varphi(n-i)q(n-i)}{2\pi \sin(2\pi i)}
$$
satisfy
\begin{equation}\label{res_cd_est}
|c_n| < C2^{-n^2}, \quad |d_n| < C2^{-n^2}
\end{equation}
with some uniform constant $C$.
Set
$$
p(z) = \sum\limits_{n \in \Z} \left( \frac{c_n}{z - n -i} + \frac{d_n}{z -n + i }\right),$$and denote by $P(t)$ the Fourier transform of $p$.

Write  $f =f_- + f_+, \,$ Sp$(f_-) \subset [-4-d,-1+d]$ and Sp$(f_+) \subset [1-d,4+d]$. Then
$$
fh = f_- h + f_+ h =  \tau_- \varphi h + \tau_+ \varphi h,
$$
where $\tau_1$ and $\tau_2$ are defined  in (\ref{tau12_def}).

Observe that the function $p(z)(\cos(2\pi z)-\cos(2\pi i))$ is entire of exponential type $2\pi$. It is also clear that the function
$\tilde\varphi:=\varphi h-p$ is entire.  Since
$$
\tilde\varphi(z)=\frac{\varphi(z)q(z)-p(z)(\cos(2\pi z)-\cos(2\pi i))}{\cos(2\pi z)-\cos(2\pi i)},
$$one may check that the exponential type of $\tilde\varphi$ is equal to $2\pi(1+d)$.
Therefore, $$\mbox{{\rm Sp}}(\tau_-\tilde\varphi)\subset \{-3,-2\}+[-1-\de,1+\de]\subset[-4-\de,-1+\de].$$ This means that $$\widehat{\tau_-\varphi h}(t)=\widehat{\tau_-(\tilde\varphi+p)}(t)=\widehat{\tau_-p}(t)=P(t-3)+cP(t-2),\quad |t|<1-\de.$$
Similarly,$$\widehat{\tau_+\varphi h}(t)=\widehat{\tau_+p}(t)=P(t+3)+cP(t+2),\quad |t|<1-\de.$$

It follows from (\ref{tau12_def}) and the proof of Lemma \ref{l5} that
$$
\widehat{\tau_-p}(t)=2\pi i\left(e^{-6\pi}+ce^{-4\pi}\right)e^{-2\pi t}\sum_n c_ne^{-2\pi i nt},\ |t|<1-\de.
$$
Similarly,
$$
\widehat{\tau_+p}(t)=-2\pi i\left(e^{-6\pi}+ce^{-4\pi}\right)e^{2\pi t}\sum_n d_ne^{-2\pi i nt},\ |t|<1-\de.
$$
By
(\ref{res_cd_est}), both series above can be extended to entire functions of $t$. Since $$\widehat{\tau_-p}(t)+\widehat{\tau_+p}(t)=0,\quad |t|<1-d,$$we conclude that \begin{equation}\label{last}\widehat{\tau_-p}(t)+\widehat{\tau_+p}(t)\equiv0.\end{equation}

 Consider  the system of functions
$$
\{w_{\pm, n}(t)\}:=\{\mathfrak{F}\left({\frac{\cos(2\pi z)-\cos(2\pi i)}{z+n\mp i}}\right), n\in\Z\},
$$where $\mathfrak{F}$ means the Fourier transform. Clearly, every function $w_{\pm,n}$ satisfies Sp$(w_{\pm,n})\subset[-1,1]$ and $w_{\pm,n}\in L^2(-1,1)$.
It is easy to check that $w_{+,n}(t)$
 is orthogonal in $L^2(-1,1)$ to every function $\exp(\pm 2\pi t +2\pi i mt)$  above except for the function $\exp(-2\pi t+2\pi i nt)$. From (\ref{last}) we deduce that  $c_n=0, n\in\Z$. Similarly, $d_n=0,n\in\Z$. Hence, $h$ does not have poles on $\Z\pm i$. Contradiction, and the lemma is proved.
\subsection{Exposed function with $\Omega(f)=\Z$}

Here we prove part (ii) of Theorem \ref{t3}.

It will be convenient for us to consider the spectra $S$ of the form
$$
S=[-5/2-\de,-1/2+\de]\cup[1/2-\de,5/2+\de],\quad 0<\de<1/2.
$$

It follows from Sec. 3.2 that there is a function $f\in$\ext\! satisfying (\ref{exp2}), (\ref{exp1}), $\Omega(f)=\Z$ and such that $f\notin$\expo. To prove part (ii) of Theorem \ref{t3}, we  will construct a function $f\in$\ext\! satisfying the conditions above and such that
 $f\in$\expo.

 The approach  is similar to the one  in Sec. 4.1, so we only present a sketch of proof. The proof consists of  a number of steps.

1. Consider the function $$\left(e^{2\pi i z}+1+e^{-2\pi i z}\right)\sin^2(\pi z)=2\left(\cos(2\pi z)+\frac{1}{2}\right)\sin^2(\pi z).$$It has only real zeros, and it has double zeros at the set of integers $\Z$ only.
Moreover, the spectrum of this function consists of points $\{\pm 2\}$ and  $\{\pm 1\}$, so that it has spectral gap $(-1,1).$

2. Let $\varphi\in PW^1_{[-1/2-\de,1/2+\de]}$ be any function with the following properties: the zeros of $\varphi$ are all real and simple, it does not vanish at the zeros of the function in Step 1, $\left\{\pm(\frac{1}{2} + \delta)\right\}$ lie in the spectrum $\varphi$, and it satisfies   \begin{equation}\label{equ1}|\varphi(x+i)|>\frac{1}{1+x^2}, \quad x\in\R,\end{equation} and \begin{equation}\label{equi2} |\varphi(n)|<2^{-n^2},\ n\in\Z.\end{equation} Such a function can be constructed as the product between a `small perturbation'    of $\sin(\pi z)$ (see Sec. 4.1)  and a  function from $PW_{[-\delta,\delta]}^1$ such that the points $\{\pm \delta\}$ lie in its spectrum.

3. Set
$$
f(z):=\left(e^{2\pi i z}+1+e^{-2\pi i z}\right)\sin^2(\pi z)\varphi(z).
$$It is easy to check that $f\in $\pw\!, the points $\{\pm (5/2+\delta)\}$ lie in its spectrum  and that   $f$ has only real zeros. By Theorem~\ref{t0}, $f\in $Ext$(PW^1_{[-5/2-\de,5/2+\de]})$.
Therefore, $f\in$\ext.

4. Using (\ref{equ1}), one may check that $f$ satisfies (\ref{exp2}). It is also easy to see that $f$ satisfies (\ref{exp1}) and that $\Omega(f)=\Z$.

5. Assume $g:=fh\in $\pw, for some nonnegative function $h$. To prove the theorem, by Lemma \ref{h_lemma}, it suffices to show that $h=const.$

Let us assume that $h=g/f\ne const.$
Similarly to above, one may check that $h$ cannot be entire function. Therefore, we may assume that $h$ is meromorphic. It is easy to see that the poles of $h$ lie on $\Z$, and that each pole is of order two. It follows that  the function $q(z):=h(z)\varphi(z)\sin^2(\pi z)$ is entire, it does not vanish at those points of $\Z$ where $h$ has poles, and one may show that
$q\in PW^1_{[-3/2-\de,3/2+\de]}.$

6. We may write
\begin{equation}\label{e0ce}
g(z)=f(z)h(z)=\left(e^{-2\pi i z}+1+e^{2\pi i z}\right)q(z).
\end{equation}Set $\Psi(t)=\hat q(t)$. Then
$$
\hat g(t)=\Psi(t-1)+\Psi(t)+\Psi(t+1), \quad t\in\R.
$$

 Recall that  $q\in PW^1_{[-3/2-\de,3/2+\de]}$. Therefore, $\Psi(t)=0, |t|\geq 3/2+\de$, and so  $$\Psi(t+k)=0,\quad  |t|\leq 1/2-\de, \ k\in\Z,\ |k|\geq 2.$$ By (\ref{e0ce}), we may write
$$
\hat g(t)=\sum_{k\in\Z}\Psi(t+k), \quad |t|<1/2-\de.
$$
Consider the Fourier series
$$
\hat g(t)\sim\sum_{k\in\Z}\alpha_ke^{2\pi i kt}, \quad |t|<1/2.
$$It is well know that the Fourier coefficients $\alpha_k$ satisfy
$$
\alpha_k=q(k),\quad k\in\Z.
$$It now follows from (\ref{equi2}) that the Fourier series above converges absolutely and  admits an analytic extension to the complex plane. However, since $g$ has spectral gap $[-1/2+\de,1/2-\de]$, i.e. $\hat g(t)=0$ for $|t|\leq 1/2-\de$, we conclude that the Fourier series vanishes identically, so that
$q(n)=0, n\in\Z$.  Contradiction.

\section{Short gap: Non-extreme function whose spectrum touches an inner border of $S$}
Recall that we consider sets of the form $S = [-\sigma, \sigma] \setminus (-\rho, \rho).$
Here we show that if the gap $(-\rho,\rho)$ is sufficiently  small then there exists a non-extreme function $f\in $\pw satisfying (\ref{c0}), $\{-\rho\}\in$Sp$(f)$ and $\La(f)=\emptyset$:
\begin{proposition}\label{t_inter_bounds}
Assume  $\sigma>7 \rho>0$. Then there exists $f\in PW_{S}^1$ which has only real zeros,  $\pm \rho\in$Sp$(f)$ and $f\not\in$\ext.
\end{proposition}

Condition $\sigma>7\rho$ above is chosen only for simplicity of the proof. It can be substituted by a less restrictive one.

\medskip
\noindent {\bf Proof.} Set
$$
f(z):= \left( \frac{\sin(\pi \rho   z)}{ z} \right)^2 \cos(4\pi\rho z).
$$
Clearly, $f \in L^1(\R)$ and $f$ has only real zeroes. Also, by the Titchmarsh convolution theorem,
$${\rm Sp}(f)\subset [-\rho,\rho]+\{-2\rho,2\rho\}= [-3\rho,-\rho]\cup[\rho,3\rho] \quad \text{ and } \quad \pm \rho \in {\rm Sp}(f),$$ so that $f\in$\pw.


To show that $f \notin$ \ext\, it suffices to provide a bounded real (on $\R$) non-constant function $h$ such that $f h \in PW^1_S$. One may choose
$$
h(z) := \cos(8\pi \rho z).
$$
Clearly,  $h$ is real, bounded and by the Titchmarsh convolution theorem,  $${\rm Sp}(fh) = \left([-3 \rho, - \rho]\cup[\rho,3\rho]\right)+\{-4\rho,4\rho\}\subset[-7\rho,-\rho]\cup[\rho,7\rho] \subset S.$$ Therefore, $fh \in PW^1_S$, which implies that $f$ is not extreme.

\section{Long gap: Proof of Theorem \ref{tt0}}

\subsection{Proof of  part (i)}

The proof will consist of several steps.

1. Recall that $f=f_++f_-,$ where Sp$(f_+)\subset [\rho,\sigma]$ and Sp$(f_-)\subset [-\sigma,-\rho]$. Denote by $\rho_\pm\geq\rho$ and $\sigma_\pm\leq\sigma$ the numbers such that Sp$(f_+)\subset [\rho_+,\sigma_+]$, $\{\rho_+\},\{\sigma_+\}\in$Sp$(f_+)$, and Sp$(f_-)\subset [-\sigma_-,-\rho_-]$, $\{-\rho_-\},\{-\sigma_-\}\in$Sp$(f_-)$.

The following claim is  easy to check
\begin{claim}\label{cl0} For every $x\in\R$ we have
\begin{equation}\label{fp}
|f_+(x+iy)|\leq Ce^{-2\pi \rho_+ y}, \ y>0,\quad |f_+(x+iy)|\leq Ce^{2\pi \sigma_+ |y|}, \ y<0,
\end{equation}where $C$ is a constant depending on $f_+$. In {\rm(\ref{fp})}, $\rho_+$ cannot be replaced by a larger constant and $\sigma_+$ cannot be replaced by smaller constant.\end{claim}

A similar claim is true for $f_-.$

Set
$$
\alpha:=\min\{\sigma-\sigma_+,\sigma-\sigma_-,\rho_+-\rho,\rho_--\rho\}.
$$Clealry $\alpha\geq0$ and
\begin{equation}\label{spe}
\alpha=\sup\{\beta\geq0: {\rm Sp}(f_\pm)+[-\beta,\beta]\subset S\}.
\end{equation}

2. Assume that $\La^\ast(f)=\emptyset$, i.e. that $f$ satisfies (\ref{la_cond}). Since $f$ is not extreme, it does not satisfy (\ref{sp_cond}), so that $\alpha>0$.  Let us use a simple argument from \cite{dLR}: Set
$$f_1(z):=(1+\cos (2\pi\alpha(z-\varphi)))f(z),\quad f_2(z):=(1-\cos (2\pi\alpha(z-\varphi)))f(z).$$Here $\varphi$ is a real number such that
$$\int_\R \cos (2\pi\alpha (x-\varphi))|f(x)|\,dx=0.$$ Such a value exists  since the integral is a continuous function of $\varphi$ which changes the sign on the interval $[0,\pi]$.
  It is clear that $f_1+f_2=2f$, both functions $f_j$ belong to  \pw and satisfy $\|f_j\|_1=1$, $\La^\ast(f_j)=\emptyset$, and that by (\ref{spe})  $f_j$ satisfy (\ref{sp_cond}). Hence, $f_j\in$\ext, $j=1,2.$

3. Assume that $\La^\ast(f)\ne\emptyset$. Consider the sequence of points $\{\lambda_j\}_{j=1}^N$, $1\leq N\leq\infty,$ such that every element $\lambda\in \La^\ast(f)\cap\{z: {\rm Im}z>0\}$, and the number of occurrences of  $\lambda$ is equal to $\min\{m(\lambda,f_-),m(\lambda,f_+)\},$
where $m(\lambda,f_-)$ is the multiplicity of zero $\lambda$ of $f_-$.

 Consider the Blaschke product (for the upper half-plane)
\begin{equation}\label{bpr}
B(z):=\prod_{j=1}^N\frac{1-z/\lambda}{1-z/\bar\lambda}.
\end{equation}

\begin{claim}\label{cl1}
Both functions $f_+B$ and $f_+/B$ belong to  $PW^1_{[\rho_+,\sigma_+]}$, and $$\{\rho_+\},\{\sigma_+\}\in{\rm Sp}(f_+B)\cap{\rm Sp}(f_+/B).$$
\end{claim}

A similar statement holds for the functions $f_-B$ and $f_-/B$.

\noindent{\bf Proof}.   Since $|B(x)|=1$ a.e. on $\R$, we see that $f_+B, f_+/B\in L^1(\R).$
Therefore, it suffices to check that both functions $f_+B$ and $f_+/B$ satisfy (\ref{fp}), where  $\rho_+$ cannot be replaced by a larger constant, and $\sigma_+$ cannot be replaced by smaller constant.

Consider the product $f_+B$. Since $|B(z)|<1,$ Im$z>0$, the first estimate in (\ref{fp}) holds true  for $f_+B$. On the other hand, one may deduce from Hayman's theorem (see \cite{levin}, Lecture 15), that the constant $\rho_+$ in that estimate cannot be replaced by a larger constant.

Since $|B(z)|>1,$ Im$z<0$, the right hand-side inequality in (\ref{fp}) cannot be true with a constant smaller than $\sigma_+$.

 If the sequence $\{\lambda_j\}_{j=1}^N$ is finite, it is clear that $f_+B$ satisfies  the right inequality in (\ref{fp}), where $\sigma_+$ cannot be replaced by a smaller constant.

Assume $N=\infty$. Set $$B_n(z):=\prod_{j=1}^n\frac{1-z/\lambda_j}{1-z/\bar\lambda_j}$$ and consider the sequence of functions
$$f_{+,n}(z):=e^{-2\pi i \sigma_+ z}f_+(z)B_n(z).$$
By the second inequality in (\ref{fp}), every $f_{+,n}$ is bounded in the lower half-plane and satisfies $|f_{+,n}(x)|\leq C$.
Therefore, by the Phragm\'{e}n–Lindel\"{o}f theorem, $|f_{+,n}(z)|\leq C$, Im$z\leq0$.
Since for every $z,$ Im$z<0$, the sequence $|B_n(z)|$ is increasing and  bounded, we conclude that $f_+B$ satisfies the right inequality in (\ref{fp}), where $\sigma_+$ cannot be replaced by a smaller constant. The claim is proved for the function $f_+B$. The rest of the proof is similar.

4. Set
$$\psi(z):=e^{2\pi i \alpha z+i\varphi }B(z)$$and$$
h(z):=\frac{1}{2}(\psi(z)+\overline{\psi(\bar z)})=\frac{1}{2}\left(B(z)e^{2\pi i \alpha z+i\varphi}+e^{-2\pi i \alpha z-i\varphi}/B(z)\right).
$$Here $\alpha$ is defined in  (\ref{spe}) and  $\varphi\in[0,\pi]$ is chosen to satisfy the condition \begin{equation}\label{hhh}\int_\R h(x)|f(x)|\,dx=0.\end{equation} Since $h$ is real on $\R$, such a number $\varphi$ exists, see  Step 2.

Clearly,
\begin{equation}\label{psi_1_est}
    |\psi(x)|=1,x\in\R, \quad \text{and} \quad |\psi(z)|\ne1, z\in \Cc\setminus\R.
\end{equation}
We  see that  $h$ is a meromorphic function which has poles at every point $\lambda\in\La^\ast(f)$, whose multiplicity is equal to $\min\{m(\lambda,f_+),m(\lambda,f_-)\}$.
We see that $hf_+$ and $hf_-$ are entire functions satisfying \begin{equation}\label{h1}|h(\lambda)f_-(\lambda)|+|h(\lambda)f_+(\lambda)|\ne0,\quad \forall \lambda\in\La^\ast(f).
\end{equation}

5. Set
$$
f_1(z):=(1+h(z))f(z),\quad f_2(z):=(1-h(z))f(z).
$$Then $f_1+f_2=2f$.

Clearly, $f_1$ and $f_2$ are entire functions integrable on $\R$. By (\ref{hhh}), $\|f_1\|_1=\|f_2\|_1=1$.
Observe that the condition $1+h(z)=0$ means that  $$1+\frac{\psi(z)+\overline{\psi(\bar z)}}{2}=\frac{\psi^2(z)+2\psi(z)+1}{2\psi(z)}=0,$$so that $\psi(z)=-1$. Similarly, the condition $1-h(z)=0$ means that $\psi(z)=1.$
 Therefore, from (\ref{psi_1_est}) and (\ref{h1}), we deduce that  $\La^\ast(f_1)=\emptyset = \La^\ast(f_2).$

By Theorem \ref{t1}, to finish the proof, we have to check that $f_1$ and $f_2$ satisfy (\ref{sp_cond}).
 Write
$$f_1(z)=(1+h(z))f(z)=(1+h(z))(f_-(z)+f_+(z)).$$

Consider the function $hf_+$.  Observe that $$
|hf_+(x+iy)|\geq \frac{1}{2}\left(e^{-2\pi\alpha y}|B(x+iy)|-e^{2\pi\alpha y}/|B(x+iy)|\right)|f_+(x+iy)|\geq$$$$\frac{1}{2}\left(e^{2\pi \alpha|y|}-1\right)|f_+(x+iy)|,\quad y<0.
$$Using  Claims \ref{cl0} and \ref{cl1}, wee see that the function $hf_+$ has exponential type $\sigma_++\alpha$ in the lower half-plane.
Similarly, $hf_+$ has the type $-\rho_++\alpha$ in the upper half-plane. Hence, $${\rm Sp}(hf_+)\subset[\rho_+-\alpha,\sigma_++\alpha]=[\rho_+,\sigma_+]+[-\alpha,\alpha],$$and $\{\rho-\alpha\},\{\sigma_++\alpha\}\in$Sp$(hf_+)$.
The same argument applied to $hf_-$ shows that Sp$(hf_-)\subset [-\sigma_-,-\rho_-]+[-\alpha,\alpha]$ and  $\{-\sigma_--\alpha\},\{-\rho_-+\alpha\}\in$Sp$(hf_-)$. By the definition of $\alpha$, we conclude that $f_1$ satisfies (\ref{sp_cond}). Similarly,
$f_2$ satisfies (\ref{sp_cond}).

\subsection{Proof of  part (ii)}

The proof consists of two steps.

1. Let $\alpha\geq0$ be the number defined in (\ref{spe}). Then both functions $e^{2\pi i \alpha z}f(z)$ and $e^{-2\pi i \alpha z}f(z)$ belong to
\pw\!, and at least one of them satisfies (\ref{sp_cond}). We may assume that $e^{2\pi i \alpha z}f(z)$ satisfies (\ref{sp_cond}).

Firstly, assume that $\La^\ast(f)=\emptyset.$
Set
$$
f_1(z):=f(z)(1-\varepsilon e^{2\pi i \alpha z}),
$$where $0<\varepsilon<1/64$ is any  number such that $f$ does not vanish on the line Im$z=-\log\varepsilon/(2\pi\alpha)$. Clearly,
$f\in $\pw\!. Since the zeros of  $1-\varepsilon\exp(2\pi i \alpha z)$ lie on the line Im$z=\log\varepsilon/(2\pi\alpha)$, we have  $\La^\ast(f_1)=\emptyset.$ It is also clear that \begin{equation}\label{eee}\|f-f_1\|_1\leq\varepsilon.\end{equation} Therefore, the function
$g_1:=f_1/\|f_1\|_1$ belongs to \ext\! and we have
$$
\|f-g_1\|_1\leq \|f-f_1\|_1+\|f_1-g_1\|_1 \leq 4 \epsilon, 
$$ which proves that   $f\in\overline{{\rm Ext}(PW_S^1)}$.

 2. Assume now that $\La^\ast(f) \neq \emptyset.$
Let $f_1,\alpha\geq 0$ and $\varepsilon>0$ be as above. Set
$$
f_2(z):=f_1(z)(1+\varepsilon B(z)).
$$
where  $B$ is the Blaschke product defined in (\ref{bpr}). Then
\begin{equation}\label{e11}\|f_2-f_1\|_1\leq\varepsilon\|f_1\|_1\leq \varepsilon(1+\varepsilon)<2\varepsilon.\end{equation} By step 1, $1-3\varepsilon\leq \|f_2\|_1\leq 1+3\varepsilon.$

\begin{claim}\label{cl2}
For every $\delta>0$ there is a positive number $\varepsilon\leq\delta$ such that the function $f_2(z)$ does not vanish on the set
$\{\lambda: \bar\lambda\in Z_c(1+\varepsilon B)\}$.
\end{claim}

The claim easily follows from the fact that  the set  $\{ B(\lambda):\lambda\in Z_c(f_2)\}$ is countable.

We assume that $\varepsilon$ is as in the claim above.

By the choice of $\alpha$ and $\varepsilon$ in step 1, $f_1$ satisfies (\ref{sp_cond}) and  $\La^\ast(f_1)=\La^\ast(f)$. Using Claim~\ref{cl1}, we see that $f_2\in$\pw.
Since $|B(z)|<1$ in the upper half-plane,  the function $1+\epsilon B(z)$ does not vanish in that plane. Moreover, $f_2$ does not vanish on the set $\La^\ast(f)\cap\{x:\mbox{Im}z<0\}$.
By Claim \ref{cl2}, $\La^\ast(f_2)=\emptyset,$ and so the function $g_2:=f_2/\|f_2\|_1$ belongs to \ext. Clearly,
$$
\|f_2-g_2\|_1\leq \|f_2\|_1|1-1/\|f_2\|_1| \leq  8 \varepsilon, 
$$Now, since $$\|f-g_2\|_1\leq \|f-f_1\|_1+\|f_1-f_2\|_1+\|f_2-g_2\|_1,$$using (\ref{eee}) and (\ref{e11}), we see that $f\in\overline{{\rm Ext}(PW_S^1)}$.

\medskip\noindent{\bf Acknowledgment}. The authors want to thank Prof. K. Dyakonov for valuable discussions about this subject. We are also grateful to Prof. M. Sodin for constructive remarks.

\noindent Alexander Ulanovskii\\ University of Stavanger, Department of Mathematics and Physics,\\
4036 Stavanger, Norway,\\
alexander.ulanovskii@uis.no

\medskip

\noindent Ilya Zlotnikov\\
University of Stavanger, Department of Mathematics and Physics,\\
4036 Stavanger, Norway,\\
ilia.k.zlotnikov@uis.no

\end{document}